\documentclass[12pt,reqno]{amsart}
\usepackage{amscd,amssymb}

\newtheorem{theorem}{Theorem}
\newtheorem{proposition}[theorem]{Proposition}

\theoremstyle{definition}

\theoremstyle{remark}

\newcommand{\Z}{{\mathbf Z}}
\newcommand{\Q}{{\mathbf Q}}
\newcommand{\R}{{\mathbf R}}
\newcommand{\C}{{\mathbf C}}
\renewcommand{\H}{{\mathcal H}}

\newcommand{\CC}{{\mathcal {C}}}

\newcommand{\lt}{{\mathfrak {\ell^2(\pi)}}}
\newcommand{\cp}{{\mathcal C_{\pi}}}
\newcommand{\ecp}{\mathcal E({\mathcal C_{\pi}})}
\newcommand{\tor}{\mathcal Tor}
\newcommand{\E}{\mathcal E}
\newcommand{\X}{\mathcal X}
\newcommand{\Y}{\mathcal Y}

\DeclareMathOperator{\im}{{im}}

\DeclareMathOperator{\ob}{{Ob}}
\DeclareMathOperator{\ho}{\hat\otimes}
\DeclareMathOperator{\spec}{Spec}

\begin{document}

\title[Zero-in-the-spectrum conjecture] %
{On the zero-in-the-spectrum conjecture}

\author[M.~Farber]{M.~Farber}
\address{Department of Mathematics, Tel Aviv University, Tel Aviv, 69978, Israel}
\email{farber@math.tau.ac.il}

\author[ S. Weinberger]{S. Weinberger}
\address{Department of Mathematics, University of Chicago, Chicago, IL}
\email{shmuel@math.uchicago.edu}

\date{October 1, 1999}

\thanks{Partially supported by the US - Israel Binational Science Foundation and by the Minkowski Center for Geometry}

\begin{abstract}
We prove that the answer to the "zero-in-the-spectrum" conjecture, in its form, suggested by J. Lott, 
is negative. Namely, we show that for any $n\ge 6$ there exists a closed $n$-dimensional smooth
manifold $M^n$,
so that zero does not belong to the spectrum of the Laplace-Beltrami operator acting on
the $L^2$ forms of all degrees on the universal covering $\tilde M$.
\end{abstract}

\maketitle

\section{\bf The Main result}

M. Gromov formulated the following conjecture (cf. \cite{G1}, p. 120;  
\cite{G2}, p.21, Problem 0.5.$F_1^{"}$, and also \cite{G2}, p. 238):

\subsection*{Conjecture A}{\it Let $M$ be a closed aspherical manifold; is it true that 
zero is
always in the spectrum of the Laplace-Beltrami operator $\Delta_p$, acting on the square integrable
$p$-forms on the universal covering $\tilde M$, for some $p$?}

If the Strong Novikov Conjecture holds 
for the fundamental group $\pi_1(M)$, then $0\in \spec(\Delta_p)$ for some $p$,
cf. \cite{L}, p. 371. Hence a conterexample to Conjecture A would be also a conterexample to the Strong
Novikov Conjecture.

G. Yu obtained in \cite{Yu1}, \cite{Yu2} results, confirming Conjecture A 
under some additional assumptions.

In 1991 J. Lott raised a more general "zero-in-the-spectrum" question:

\subsection*{Conjecture B}{\it Is it true, that for any complete Riemannian manifold $M$ zero is
always in the spectrum of the Laplace-Beltrami operator $\Delta_p$, acting on the square integrable
$p$-forms on $M$, for some $p$?}

We refer to the survey articles \cite{L} and \cite{Lu}.

J. Lott showed in \cite{L} that Conjecture B is true for manifolds of low dimension and also
for some classes of higher dimensional manifolds.

In this article we give negative answers to Conjecture B and also to a version of Conjecture A
where one removes the assumption of asphericity of $M$.
We prove the following Theorem.

\begin{theorem}\label{theorem1}  For any $n\ge 6$
there exists a closed $n$-dimensional smooth manifold $M$, 
so that for any $p=0, 1, \dots, n$
the zero does not belong to the spectrum of the Laplace-Beltrami operator
\[\Delta_p:\, \,  \Lambda^p_{(2)}(\tilde M) \to \Lambda^p_{(2)}(\tilde M),\]
acting on the space of $L^2$-forms $\Lambda^p_{(2)}(\tilde M)$ 
on the universal covering $\tilde M$
of $M$.\end{theorem}

Our proof of Theorem 1 will be based on the fact that it
can be restated in an equivalent form using the notion of extended $L^2$-homology, introduced in \cite{F1}: 

\begin{theorem}\label{theorem2} For any $n\ge 6$
there exists a closed orientable $n$-dimensional smooth manifold $M$, 
so that extended $L^2$ homology $\H_p(M;\lt)=0$ vanishes for all $p$. Here 
$\pi$ denotes the fundamental group $\pi=\pi_1(M)$, and $\lt$ denotes the $L^2$-completion
of the group ring $\C[\pi]$.\end{theorem}

Equivalence between Theorem \ref{theorem1} and Theorem \ref{theorem2} can be 
established as follows.
Zero not in the spectrum of the Laplacian 
$\Delta_p:\, \,  \Lambda^p_{(2)}(\tilde M) \to \Lambda^p_{(2)}(\tilde M)$ for all $p$ 
is equivalent to vanishing of the extended $L^2$ cohomology $\H^\ast(M;\lt)$, cf. \cite{F1},
according to the De Rham Theorem for extended cohomology, cf. section 7 of \cite{F2} and also \cite{S}.
Vanishing of the extended $L^2$-cohomology is equivalent to vanishing of the extended 
$L^2$-homology 
$\H_\ast(M;\lt)$, because of the Poincar\'e
duality, cf. \cite{F1}, Theorem 6.7.

The proof of Theorem \ref{theorem2} is based on the following Theorem:
\begin{theorem}\label{thm3}
There exists a finite 3-dimensional polyhedron $Y$ with fundamental group 
$\pi_1(Y)=\pi =F\times F\times F$, where $F$ denotes a free group with two generators, such that
the extended $L^2$-homology $\H_p(Y;\ell^2(\pi))=0$ vanishes for all $p=0, 1, \dots$.
\end{theorem}

The strategy of our proof of Theorems \ref{theorem2} and \ref{thm3}
is similar to the method used by M.A. Kervaire \cite{K}, 
who constructed smooth homology spheres with prescribed fundamental groups. Our proof
uses $L^2$-analog of the Hopf exact sequence.

The authors are thankful to B. Eckmann and A. Connes for helpful conversations.

\section{\bf Proofs of Theorems 2 and 3}

{\bf A.} Let $\pi$ be a discrete group given by a finite presentation
\[\pi = |x_1, x_2, \dots, x_n: r_1=1, r_2=1, \dots, r_m=1|\]
by generators and relations. We will assume that:
\begin{enumerate}
\item[(a)] The extended $L^2$-homology of $\pi$ in dimensions $0, 1$ and $2$ vanishes, i.e.
$$\H_0(\pi;\lt) = \H_1(\pi;\lt) = \H_2(\pi;\lt) =0.$$
\item[(b)] Let $X$ be a finite cell complex with $\pi_1(X)=\pi$ having one zero-dimensional cell,
$n$ cells of dimension 1 and $m$ cells of dimension two, constructed in the usual way out of 
the given presentation of $\pi$. 
Then the second homotopy group $\pi_2(X)$ of $X$, 
viewed as a $\Z[\pi]$-module, is free and finitely generated. 
\end{enumerate}

Our purpose is to show that there exists a 3-dimensional cell complex
$Y$, obtained from $X$ by first taking a bouquet with finitely many copies of $S^2$ and then
adding a finite number of 3-dimensional cells, so that 
\begin{eqnarray}
\H_i(Y;\lt) = 0\quad\text{for any}\quad i=0, 1, \dots.
\end{eqnarray}

{\bf B. $L^2$-Hopf exact sequence.}
First we will calculate the extended $L^2$ homology of $X$ using the spectral sequence constructed
in Theorem 9.7 of \cite{F1}.
We will work in the von Neumann category $\cp$
of Hilbert representations of $\pi$, cf. \cite{F2}, \S 2, example 5. We will denote by $\ecp$ the corresponding extended abelian category, cf. \cite{F2}, \S 1. Let $\tilde X$ be the universal covering
of $X$. 
We will use the functors
\[\tor_p^{\pi}(\lt, H_q(\tilde X;\C))\] 
with values in the extended abelian category $\ecp$, which are defined in \cite{F1}, page 660 under the 
assumption that the homology modules $H_q(\tilde X;\C)$ of the universal covering admit finite
free resolutions. In our case only two of these homology modules can be nonzero (for $q=0$ and 
$q=2$), and (since $H_2(\tilde X;\C) = \C\otimes \pi_2(X)$) our assumption (b) guarantees this finiteness condition for $q=2$. 
The functor $\tor_0^\pi(\lt, H_q(\tilde X;\C))$ can be denoted by
\begin{eqnarray}
\lt\tilde\otimes_\pi H_q(\tilde X;\C).
\end{eqnarray}
It is an analog of the tensor product, cf. \cite{F2}, \S 6.
Note that in general it takes values in the extended category
$\ecp$, i.e. it may have a nontrivial torsion part.

By Theorem 9.7 of \cite{F1}, there exists
a spectral sequence in the abelian category $\ecp$ with the following properties:
\begin{itemize}
\item the initial term $E^2_{p,q}$ equals $\tor^\pi_p(\lt, H_q(\tilde X;\C))$. 
\item The spectral sequence converges to the extended homology $\H_{p+q}(X;\lt)$.
\end{itemize}

For $q=0$ we have $H_0(\tilde X;\C)=\C$, and $\tor^\pi_q(\lt, \C)$ can be also understood as
the extended $L^2$ homology of the Eilenberg - MacLane space $K(\pi, 1)$.
We will use notation
\begin{eqnarray}
\tor^\pi_q(\lt, \C) = \H_q(\pi; \lt).
\end{eqnarray}
It is an analog of the group homology.

Since $X$ is two-dimensional, 
the spectral sequence contains only two rows ($q=0$ and $q=2$) and may have only one nontrivial differential. Hence
we obtain the following isomorphisms:
\begin{eqnarray}
\H_0(X;\lt) \simeq \H_0(\pi;\lt)\quad \text{and}\quad \H_1(X;\lt) \simeq \H_1(\pi;\lt).\label{hur0}
\end{eqnarray}
These are Hurewicz type isomorphisms. 
The first nontrivial differential of the $E^2$-term is $d_2: E^2_{3,0}\to E^2_{0,2}$.
Here $E^2_{3,0} = \H_3(\pi;\lt)$ and $E^2_{0,2}= \lt\tilde\otimes_\pi H_2(\tilde X;\C)$.
Using the Hurewicz isomorphism $H_2(\tilde X)\simeq \pi_2(\tilde X)\simeq \pi_2(X)$,
we may write
\begin{eqnarray*}
E^2_{0,2}= \lt\tilde\otimes_\pi \pi_2(X)
\end{eqnarray*}
and the above differential is 
\begin{eqnarray}
d_2: \H_3(\pi;\lt) \to \lt\tilde\otimes_\pi \pi_2(X).
\end{eqnarray}
Note also that this differential must be a monomorphism (viewed as a morphism of
the abelian category $\ecp$), since $\H_3(X;\lt)=0$ 
(recall that $X$ is two-dimensional).
The spectral sequence above yields the following exact sequence
\begin{eqnarray}\quad\quad  0\to \H_3(\pi;\lt)\stackrel{d_2}{\to}\lt\tilde\otimes_\pi \pi_2(X) \stackrel{h}{\to} 
\H_2(X;\lt) \to \H_2(\pi,\lt) \to 0.\label{hopf}
\end{eqnarray}
It is an $L^2$ analog of the {\it Hopf's exact sequence}. 

We conclude (using (\ref{hur0}) and our assumptions (a))
that 
$$\H_0(X;\lt) = \H_1(X;\lt) =0$$
 and $\H_2(X;\lt)$ can be found from the exact sequence
\begin{eqnarray}\quad\quad  0\to \H_3(\pi;\lt)\to\lt\tilde\otimes_\pi \pi_2(X) \stackrel{h}{\to} 
\H_2(X;\lt) \to 0.\label{hopf1}
\end{eqnarray}

{\bf C.} We will now specialize our discussion to the following group 
$$\pi=F\times F\times F,$$ 
where 
$F$ is a free group with two generators. We will denote the free generators of the factor number $r$
(where $r=1, 2, 3$) 
by $a^r_1, a^r_2$. 
We will fix the presentation of $\pi$ given by 6 generators
$a^1_1, a^1_2, a^2_1, a^2_2, a^3_1, a^3_2$ and the following 12 relations
\[(a^k_i,a^l_j)=1,\quad \text{for}\quad k\ne l, \quad k,l\in\{1, 2, 3\}, \quad i,j\in \{1,2\},\]
where $(v,w)=vwv^{-1}w^{-1}$ denotes the commutator. 

$\pi$ satisfies condition (a) above, as follows from the Kunneth theorem for the extended
$L^2$-cohomology, cf. Appendix, Theorem \ref{thm5}. Here we use that $\H_j(F;\ell^2(F))$
is nonzero only for $j=1$ and
has no torsion; hence the terms containing the periodic product in formula (\ref{kunneth2}),  
vanish; cf. Proposition \ref{prop3}, statement (b).

Let us show that this group $\pi$, together with its specified 
presentation, satisfies condition (b).
The two-dimensional complex $X$ constructed out of this presentation 
will have one zero-dimensional cell $e^0$, six 1-dimensional cells $e^1_i, e^2_i, e^3_i$ 
and 12 two-dimensional cells
$e^{12}_{ij}, e^{13}_{ij}, e^{23}_{ij}$. Here $e^k_i$ denotes the 1-cell corresponding to the generator $a^k_i$ and $e^{kl}_{ij}$ denotes the 2-cell corresponding to the relation $(a^k_i,a^l_j)=1$.

Let $0\to C_2\to C_1\to C_0\to 0$ denote the chain complex of the universal covering 
$\tilde X$. The boundary homomorphism acts as follows
$$
\begin{array}{l}
\partial e^k_i =(a^k_i-1)e^0,\\
\partial e^{kl}_{ij} =(a^k_i-1)e^l_j - (a^l_j-1)e^k_i.\\
\end{array}
$$

Using the Hurewicz isomorphisms $\pi_2(X)=\pi_2(\tilde X)=H_2(\tilde X)$,
we may compute the group $\pi_2(X)$, viewed as a $\Z[\pi]$-module, as the kernel of 
$\partial: C_2\to C_1$. Let 
\[x\in C_2, \quad 
x=\sum_{ij}\lambda^{12}_{ij}e^{12}_{ij} + \sum_{ij}\lambda^{13}_{ij}e^{13}_{ij} + \sum_{ij}\lambda^{23}_{ij}e^{23}_{ij},\quad \lambda^{kl}_{ij}\in \Z[\pi],\]
be an element with $\partial x=0$. Then the following equations hold
$$
\begin{array}{l}
\sum_i \lambda^{12}_{ij}(a^1_i-1) =\sum_i \lambda^{23}_{ji}(a^3_i-1),\\ \\
\sum_j \lambda^{12}_{ij}(a^2_j-1) + \sum_j \lambda^{13}_{ij}(a^3_j-1)=0,\\ \\
\sum_i \lambda^{13}_{ij}(a^1_i-1) + \sum_i \lambda^{23}_{ij}(a^2_i-1) =0.\\
\end{array}
$$
Hence we may write
$$
\begin{array}{l}
\lambda^{12}_{ij} = \sum_k \mu^{12}_{ijk}(a^3_k-1), \quad \mu^{12}_{ijk}\in \Z[\pi],\\ \\
\lambda^{23}_{ij} = \sum_k \mu^{23}_{ijk}(a^1_k-1), \quad \mu^{23}_{ijk}\in \Z[\pi],\\ \\
\lambda^{13}_{ij} = \sum_k \mu^{13}_{ijk}(a^2_k-1), \quad \mu^{13}_{ijk}\in \Z[\pi].\\ \\
\end{array}
$$
Therefore we obtain
\begin{eqnarray}
\mu^{12}_{ijk} = \mu_{jki}^{23} = - \mu_{ikj}^{13}.\label{equ1}
\end{eqnarray}Conversely, any system $\mu_{ijk}^{rs}\in \Z[\pi]$ satisfying (\ref{equ1}) 
determines a cycle 
$x\in C_2$, $\partial x=0$.
This proves that {\it $\pi_2(X)$ is a free $\Z[\pi]$-module of rank 8 with the basis}
\begin{eqnarray}
x_{ijk}=(a^1_i-1)e^{23}_{jk} - (a^2_j-1)e^{13}_{ik}+(a^3_k-1)e^{12}_{ij}, \quad i,j,k\in \{1, 2\}.
\label{equ2}
\end{eqnarray}

Note that the Eilenberg-MacLane complex $K=K(\pi,1)$ is $B\times B\times B$, where $B$ is the bouquet of two circles; $K$ is obtained from $X$ by adding 8 
three-dimensional cells $e_{ijk}$, where $i,j,k\in \{1,2\}$, 
which correspond to different triple products of 1-dimensional cells of $B$. 
It is easy to see that
the boundary of $e_{ijk}$ is given by
\[\partial e_{ijk} = x_{ijk}\in \pi_2(X).\]
The chain complex of the universal covering $\tilde K$ is 
$0\to C_3\to C_2\to C_1\to C_0\to 0$, where $C_3$ is the free $\Z[\pi]$-module
generated by the cells $e_{ijk}$ and the rest is the same as the chain complex of $\tilde X$.

For a discrete group $\pi$ we will denote by $C^\ast_\R(\pi)\subset C^\ast_r(\pi)$ 
the real part of the
reduced $C^*$-algebra, i.e. the norm closure of the real group ring $\R[\pi]\subset \C[\pi]$.

{\bf D. Proposition.} {\it
Let $F$ be the free group with generators $a_1, a_2$. 
Then there exist $u_1, u_2\in C^\ast_\R(F)\subset C^\ast_r(F)$ such that 
\begin{enumerate}
\item[(i)] $u_1(a_1-1) +u_2(a_2-1)=0,$
\item[(ii)] for any pair $v_1, v_2\in \ell^2(\pi)$ with 
\begin{eqnarray}
v_1(a_1-1) +v_2(a_2-1)=0\label{equa}
\end{eqnarray} 
there exists a unique
$w\in \ell^2(\pi)$ such that 
\[v_1=wu_1, \quad v_2=wu_2.\]
\end{enumerate}
Here we consider $F$ as a subgroup of $\pi= F\times F\times F$ identifying it with one of the
factors.The reduced $C^\ast$-algebra $C^\ast_r(F)\subset C^\ast_r(\pi)$ acts in the usual way on $\ell^2(\pi)$.}
\begin{proof} For convenience, we will assume in the proof that $F$ is the third factor in $\pi$.
Consider the standard complex
\begin{eqnarray}
\ell^2(F)\oplus \ell^2(F)\stackrel{d}\to \ell^2(F), \quad (v_1, v_2)\mapsto v_1(a_1-1)+v_2(a_2-1).\label{equa3}
\end{eqnarray}
calculating the extended $L^2$ homology of the bouquet $S^1\vee S^1$
of two circles with coefficients in $\ell^2(F)$.
Since $F$ is not amenable, we know from Brooks' theorem that $d$ is an epimorphism, i.e. 
$\H_0(S^1\vee S^1;\ell^2(F))=0$. The Euler characteristic calculation shows that 
$\ker d =\H_1(S^1\vee S^1;\ell^2(F))$ is one dimensional, i.e. it is isomorphic to $\ell^2(F)$.
Here we use the fact that the von Neumann algebra $\mathcal N(F)$ is a factor. 

Let $P: \ell^2(F)\oplus \ell^2(F)\to \ell^2(F)\oplus \ell^2(F)$ be the orthogonal projection onto 
$\ker d$. We claim  that the element $P(1,0)$ belongs to 
\begin{eqnarray}
C^\ast_\R(F)\oplus C^\ast_\R(F)\subset 
C^\ast_r(F)\oplus C^\ast_r(F)\subset \ell^2(F)\oplus \ell^2(F).\label{reduced}
\end{eqnarray}
Let $d^\ast$ be the adjoint of $d$. Then $\ker d = \ker (d^\ast d)$. Moreover, the image of $d^\ast d$
is closed and thus zero is an isolated point in the spectrum of $d^\ast d$. 
Hence we may use the 
holomorphic functional calculus (Cauchy's formula)
in order to express the projector $P$ as
\begin{eqnarray*}
P = \frac{1}{2\pi i} \int_\Gamma (z-d^\ast d)^{-1}dz,
\end{eqnarray*}
where $\Gamma$ is a small circle around the origin. This explains that $P(v_1,v_2)$
belongs to $C^\ast_r(F)\oplus C^\ast_r(F)$ (cf. (\ref{reduced})), 
assuming that $v_1,v_2$ lie in the reduced $C^\ast$-algebra
$C^\ast_r(F)$. Moreover, since the operator $d^\ast d$ is real, we obtain that $P(v_1,v_2)\in \, C^\ast_\R(F)\oplus C^\ast_\R(F)$, 
for $v_1,v_2\in C^\ast_\R(F)$. 

We will set now 
$$(u_1,u_2) =P(1,0).$$
Then (i) is clearly satisfied.

We want to show that the restriction of $P$ on the first summand $\ell^2(F)$ in (\ref{equa3})
gives an isomorphism
$P: \ell^2(F)\to \ker d$. Since both $\ker d$ and $\ell^2(F)$ have von Neumann dimension one,
and the spectrum of
$P$ contains only $0$ and $1$, we conclude that it is enough to show that $P(v,0)=0$ for 
$v\in \ell^2(F)$ implies $v=0$. 
If $P(v,0)=0$ i.e. $(v,0)\in ({\ker d})^{\perp}$ then 
$\langle v, \ker d\rangle =0$, i.e. $v$ is orthogonal
to the projection of $\ker d$ on the first summand $\ell^2(F)$. 
From this we will obtain that necessarily
$v=0$  
if we will show that the projection of $\ker d$ on the first summand is dense. 

Let $f_i: \ell^2(F)\to \ell^2(F)$ be operator $x\mapsto x(a_i-1)$, where $i=1, 2$. It is clear 
that $f_1$ and $f_2$ are injective and hence their images are dense. We claim that $f_1^{-1}(\im f_2)$
is dense. If not, let $H$ denote the orthogonal complement to the closure of  
$f_1^{-1}(\im f_2)$. Then we may apply Proposition 2.4 from \cite{F2}; it implies that
$H$ must intersect $\im f_2$, which is impossible. 
Hence it follows that the projection of $\ker d$ on the first summand $\ell^2(F)$ 
(which coincides with $f_1^{-1}(\im f_2)$) is dense.

As a result we obtain from the above arguments that for any pair $(v_1, v_2)\in \ker d$ 
(i.e. which is a solution of (\ref{equa})) there exists $w\in \ell^2(F)$, so that $P(w,0)=(v_1,v_2)$, i.e. $v_1=wu_1$ and $v_2=wu_2$.
This is in fact a part of our statement (ii).
 
In order to prove (ii) in full generality, observe that
\begin{eqnarray}
\ell^2(\pi) = \ell^2(F)\hat \otimes \ell^2(F)\hat\otimes \ell^2(F),
\end{eqnarray} 
(cf. appendix) and thus (using the Kunneth theorem for extended $L^2$ homology, cf. Theorem \ref{thm4})
we find that the kernel of the operator
$$d: \ell^2(\pi)\oplus \ell^2(\pi)\to \ell^2(\pi), \quad (v_1, v_2)\mapsto v_1(a_1-1)+v_2(a_2-1),$$
equals $\ell^2(F)\hat\otimes \ell^2(F)\ho\H_1(S^1\vee S^1; \ell^2(F))$. 
(ii) now follows. \end{proof}

{\bf E.} Now we describe the kernel of the Hurewicz homomorphism 
$$h: \lt\tilde\otimes_\pi \pi_2(X)\to \H_2(X;\lt).$$
Let $u^s_i\in C^\ast_r(\pi)$, where $s=1, 2, 3$ and $i=1, 2$,
denote the element given by Proposition {\bf D} applied to the factor $F\subset \pi$ number 
$s=1, 2, 3$. Here we consider $C^\ast_r(F)$ as being canonically embedded into the von Neumann algebra $\mathcal N(\pi)$. 

We claim that {\it the kernel of the Hurewicz homomorphism $h$ is generated
by the element 
\begin{eqnarray}
y = \sum_{ijk}u^1_iu^2_ju^3_kx_{ijk}\, \in \, C^\ast_\R(\pi)\tilde 
\otimes_\pi \pi_2(X).\label{element}
\end{eqnarray}
More precisely, our statement is that any element 
$x \in\,  \lt\tilde\otimes_\pi \pi_2(X)$ with $h(x)=0$ has the form $x=\mu y$ for some} $\mu\in \ell^2(\pi)$. 

Note that the product $\mu y$ has sense because the coefficients of $y$ in the basis
$x_{ijk}$ belong to $C^\ast_\R(\pi)\subset C^\ast_r(\pi)$. 

First, it is easy to check (using (\ref{equ2})) that $h(y)=0$. 

Let 
$$x=\sum_{ijk} \mu_{ijk} x_{ijk}\, \in\,  \lt\tilde\otimes_\pi \pi_2(X), \quad h(x)=0,$$ 
be an arbitrary element of $\ker h$, where $\mu_{ijk}\in \ell^2(\pi)$. Using (\ref{equ2}), we 
obtain (equating to zero the coefficients of the cells $e^{23}_{jk}$)
that for any pair of indices $j, \, k$ holds
$$\sum_{i=1}^2 \mu_{ijk}(a_i^1-1) =0.$$
Hence, applying Proposition {\bf D}, we conclude that there exist $\mu_{jk}\in \ell^2(\pi)$
such that 
\begin{eqnarray}
\mu_{ijk}=\mu_{jk}u^1_i.\label{equ3}
\end{eqnarray}

We write again $h(x)=0$, equating to zero the coefficients of the cells $e^{13}_{ik}$ and using (\ref{equ3}). We obtain that for any pair of indices $i, \, k$ holds
\begin{eqnarray}
\sum_j \mu_{jk}u_i^1(a_j^2-1)\,  =\,  \left[\sum_j \mu_{jk}(a_j^2-1)\right]u_i^1 \, =\, 0.\label{equ4}
\end{eqnarray}
Note that $wu_1^s=0$ for $w\in \ell^2(\pi)$ implies $wu_2^s=0$ (using (\ref{equa})) and 
from the uniqueness statement in Proposition {\bf D}, (ii), we obtain that $w=0$.
Therefore (\ref{equ4}) 
implies
\begin{eqnarray*}
\sum_j \mu_{jk}(a_j^2-1)=0
\end{eqnarray*}
and hence using Proposition {\bf D}, 
\begin{eqnarray*}
\mu_{jk}=\mu_k u_j^2, \quad\rm{where}\quad \mu_k\in \ell^2(\pi).
\end{eqnarray*}
Substitute again $\mu_{ijk}=\mu_k u^1_i u^2_j$ into $h(x)=0$ and 
equating to zero the coefficients of the cells $e^{13}_{ik}$ we obtain
\begin{eqnarray}
\left[\sum_k \mu_{k}(a_k^3-1)\right]u^1_i u_j^2=0, 
\quad{\rm and\,  hence}\quad \sum_k \mu_{k}(a_k^3-1)=0.
\end{eqnarray}
Using Proposition {\bf D} as above we finally obtain 
\begin{eqnarray*}
\mu_{k}=\mu u_k^3, \quad\rm{where}\quad \mu\in \ell^2(\pi).
\end{eqnarray*}
Therefore, we find that $\mu_{ijk}=\mu u_i^1u_j^2u_k^3$ and $x=\mu y$. \qed

{\bf F.} Our goal is to show that {\it one may add 8 cells of dimension 3 to the bouquet $X\vee S^2$
such that the obtained 3-dimensional complex $Y$ will have all trivial extended $L^2$ homology}
$$\H_j(Y;\lt) =0, \quad j=0, 1, \dots$$.

For the proof, let's examine again the exact sequence (\ref{hopf1}):
\begin{eqnarray}
0\to \H_3(\pi;\lt)\stackrel{\phi}\to\lt\tilde\otimes_\pi \pi_2(X) \stackrel{h}{\to} 
\H_2(X;\lt) \to 0.\label{hopf2}
\end{eqnarray}
As we know, $\phi$ maps the generator $y$ of $\H_3(\pi;\lt)$ according to formula (\ref{element}),
i.e. $\phi$ is given by a matrix with entries in $C^\ast_\R(\pi)\subset C^\ast_r(\pi)$. 
Let $$Q: \lt\tilde\otimes_\pi \pi_2(X)\to \lt\tilde\otimes_\pi \pi_2(X)$$
denote the orthogonal projection onto $(\im \phi)^{\perp}$, the orthogonal complement of the
image of $\phi$. Since $X$ is two-dimensional, $\H_2(X;\lt)$
in has no torsion and therefore $\im \phi$ is closed. 
Note that $(\im \phi)^{\perp}$ coincides with $\ker (\phi\phi^\ast)$.
Since the image of $\phi\phi^\ast$ is closed we conclude that zero is an isolated point in the
spectrum of $\phi\phi^\ast$ and hence we may write
$$Q = \frac{1}{2\pi i}\int_\Gamma (z-\phi\phi^\ast)^{-1}dz,$$
where $\Gamma$ is a small circle round zero. Therefore, 
{\it in the basis $x_{ijk}$ the projector $Q$ is given a $(8\times 8)$-matrix with entries in $C^\ast_\R(\pi)$. }

The projective $C^\ast_\R(\pi)$-module determined by $Q$ is stably free; we know that adding a
free one-dimensional module (generated by $y$)
makes it free. Therefore we may consider the bouquet
$X_1 = X\vee S^2$ so that $\H_2(X_1;\lt)=\H_2(X;\lt)\oplus \lt$ and $\pi_2(X_1) = \pi_2(X)\oplus
\Z[\pi]$. Thus, the exact sequence (\ref{hopf2}) for $X_1$ 
\begin{eqnarray}
0\to \H_3(\pi;\lt)\stackrel{\psi}\to\lt\tilde\otimes_\pi \pi_2(X_1) \stackrel{h}{\to} 
\H_2(X_1;\lt) \to 0.\label{hopf3}
\end{eqnarray}
will have the following property: {\it the orthogonal projection 
$$Q_1: \lt\tilde\otimes_\pi \pi_2(X_1)\to \lt\tilde\otimes_\pi \pi_2(X_1)$$
onto $(\im \psi)^{\perp}$ is given by a $(9\times 9)$-matrix with entries in $C^\ast_\R(\pi)$
which determines a free $C^\ast_\R(\pi)$-module of rank 8}.

 We may reformulate the last statement as follows: there exists a $\Z[\pi]$-homomor\-phism
\begin{eqnarray}
\gamma: (\Z[\pi])^8\to C_\R^\ast(\pi)\tilde\otimes_{\pi}\pi_2(X_1)\label{homo}
\end{eqnarray}
such that the following composite
\begin{eqnarray}
\lt\tilde\otimes_{\pi}(\Z[\pi])^8 \stackrel{1\otimes \gamma}\to \lt \tilde\otimes_{\pi}\pi_2(X_1)\stackrel
{h}\to \H_2(X_1;\lt)\label{compos}
\end{eqnarray}is an isomorphism. 
Now we will use the fact that the rational group ring $\Q[\pi]$ is dense in 
$C^\ast_\R(\pi)$ with respect to the operator norm topology. Hence we may approximate
$\gamma$ by a $\Z[\pi]$-homomorphism
\begin{eqnarray*}
\gamma_1: (\Z[\pi])^8\to \Q[\pi]\otimes_{\pi}\pi_2(X_1)\label{homo1}
\end{eqnarray*}
so that the similar composition (\ref{compos}) is an isomorphism. Finally, we may multiply 
$\gamma_1$ by a large integer $N$ to obtain a $\Z[\pi]$-homomorphism
\begin{eqnarray*}
\gamma_2: (\Z[\pi])^8\to \Z[\pi]\otimes_{\pi}\pi_2(X_1)=\pi_2(X_1)\label{homo2}
\end{eqnarray*}
such that the composition 
\begin{eqnarray}
(\lt)^8 = \lt\tilde\otimes_{\pi}(\Z[\pi])^8 \stackrel{1\otimes \gamma_2}\to \lt \tilde\otimes_{\pi}\pi_2(X_1)\stackrel
{h}\to \H_2(X_1;\lt)\label{compos1}
\end{eqnarray}
is an isomorphism. 

Let $z_1, \dots, z_8\in \pi_2(X_1)$ be images of a free basis of $(\Z[\pi])^8$ under $\gamma_2$.
Realize each $z_j$ by a continuous map $f_j: S^2\to X_1$, where $j=1, \dots, 8$, and let 
$$Y = X_1\cup e_1^3\cup\dots \cup e^3_8$$
be obtained from $X_1$ by glueing $8$ three-dimensional cells to $X_1$ along $f_1, \dots, f_8$.
We claim that
\begin{equation}
\H_j(Y;\lt) =0\quad{\rm for \, \, all}\quad j=0, 1, .\dots\label{zero}
\end{equation}
In order to show this, we note that $\H_j(Y,X;\lt)$ vanishes for all $j\ne 3$ and the 3-dimensional extended $L^2$ homology $\H_3(Y,X;\lt)$
equals $(\ell^2(\pi))^8$. The boundary homomorphism
$\partial: \H_3(Y,X;\lt)\to \H_2(X;\lt)$ is an isomorphism since it coincides with (\ref{compos1}).
Hence (\ref{zero}) follows from the homological exact sequence of the pair $(Y,X)$.
This completes the proof of Theorem \ref{thm3}.

{\bf G.} Now we may complete the proof of Theorem 2. We have 
constructed above a finite 3-dimensional 
polyhedron $Y$. For any $n\ge 6$ we may embed $Y$ into $\R^{n+1}$ as a subpolyhedron.
Let $N\subset \R^{n+1}$ be the regular neighborhood of $Y\subset \R^{n+1}$. We will define
$M$ as the boundary of $N$, i.e. $M=\partial N$.

First note that the inclusion $M\to N$ induces an isomorphism of the fundamental groups and thus
$\pi_1(M) = \pi = F\times F\times F,$
where $F$ is a free group in two generators.
We want to show that
\begin{equation}
\H_j(M;\lt) =0,\quad{\rm for\, all}\quad  j=0, 1, \dots\label{zero1}
\end{equation}
In the exact homological sequence 
$$\dots \to \H_j(M;\lt)\to \H_j(N;\lt) \to \H_j(N,M;\lt) \to \dots$$
we have $\H_j(N;\lt)=0$. Also, $\H_j(N,M;\lt)\simeq \H^{n+1-j}(N;\lt)$ 
by the Poincar\'e duality 
(cf. \cite{F1})
and $\H^{n+1-j}(N;\lt)=0$ because of (\ref{zero})
using the Universal Co\-efficients Theorem (cf. \cite{F1}). 
Hence,  (\ref{zero1}) follows.
\qed

\section*{\bf Appendix: Kunneth theorem for extended $L^2$ cohomology}

\subsection*{1.} A {\it Hilbert category} $\CC$ is defined as 
an additive subcategory of the category of Hilbert spaces
and bounded linear maps, 
such for any morphism $f: H\to H'$ of $\CC$ the inclusion $\ker(f)\subset H$
belongs to $\CC$ and also the adjoint map $f^\ast: H'\to H$ belongs to $\CC$, cf. \cite{F1}.
It is shown in \cite{F1} that any Hilbert category can be canonical embedding into an abelian
category $\E(\CC)$, called {\it the extended abelian category}. 

Let $\CC$, $\CC'$ and $\CC''$ be three Hilbert categories and let 
\begin{equation}
\ho: \CC\times \CC'\, \to \, \CC''\label{tensor}
\end{equation}
be a covariant functor of two variables ({\it the "tensor product"}) such that 
\begin{enumerate}
\item[(a)] for $H\in \ob(\CC)$ and $H'\in \ob(\CC')$ the image $H\ho H'$ has as 
the underlying Hilbert 
space the tensor product of Hilbert spaces $H$ and $H'$;
\item[(b)] if $f: H\to H_1$ is a morphism of $\CC$ and $f': H'\to H_1'$ is a morphism of $\CC'$ then
$f\ho f': H\ho H'\to H_1\ho H'_1$ is the tensor product of bounded linear maps $f$ and $f'$.
\end{enumerate}

Recall that the tensor product if Hilbert spaces $H\ho H'$ is defined as the Hilbert space completion
of the algebraic tensor product $H\otimes H'$ with respect to the following scalar product
$\langle v\otimes w, v'\otimes w'\rangle = \langle v,v'\rangle\cdot \langle w,w'\rangle.$

Suppose that $(C, d)$ and $(C',d)$ are chain complexes in $\CC$ and
$\CC'$ correspondingly. 
We assume that all chain complexes are graded by non-negative integers
and have a finite length.
Their tensor product $(C,d)\ho (C',d)$ (defined in the usual way)
is a chain complex in $\CC''$. $(C,d)\ho (C',d)$ is a projective chain complex 
in the abelian category $\E(\CC'')$ and
its extended 
homology $\H_\ast(C\ho C')$ is an object of the extended category $\E(\CC'')$.
Our purpose is to express the extended homology of $(C,d)\ho (C',d)$ in terms of the extended
homology $\H_\ast(C)$ of $(C,d)$ and $\H_\ast(C')$ of $(C',d)$.

{\bf 2. Example.} Suppose that $G$ and $H$ are discrete groups. Let $\CC_G$ 
denote the category of Hilbert
representations of $G$. Recall, that an object of $\CC$ is a Hilbert space with a unitary $G$-action 
which can be continuously and 
$G$-equivariantly imbedded  into a finite direct sum $\ell^2(G)\oplus\dots \oplus \ell^2(G)$; morphisms
of $\CC$ are bounded linear maps commuting with the $G$-action. 
Then we have the tensor product functor
\begin{equation}
\ho: \CC_G\times \CC_H\, \to \, \CC_{G\times H}\label{tensor1}
\end{equation}
which is of a primary interest for us.

\subsection*{3. Tensor and periodic products} 
Given a tensor product (\ref{tensor}), it defines
two bifunctors $\E(\CC)\times \E(\CC)\to \E(\CC)$, which we now describe. 
Let $\X=(\alpha:A'\to A)\in \ob(\E(\CC))$ and $\Y=(\beta:B'\to B)\in \ob(\E(\CC'))$
be two objects with $\alpha$ and $\beta$ injective. 
Consider the following chain complex in $\CC''$
\begin{equation}
0\to A'\ho B'\stackrel{\left(
\begin{array}{c}
-1\ho\beta \\
\alpha\ho 1
\end{array}\right)
}\longrightarrow (A'\ho B) \oplus (A\ho B')
\stackrel{\displaystyle{(\alpha\ho 1, 1\ho\beta)}}\longrightarrow A\otimes B\to 0. \label{chain}
\end{equation}In other words, we view the objects $\X$ and $\Y$ as chain complexes of 
length 1 and then (\ref{chain}) is the tensor product of these chain complexes. 
The extended 
homology of (\ref{chain}) in dimension 0 will be called {\it the tensor product of
$\X$ and $\Y$}: 
\begin{eqnarray}
\X\ho\Y = ((\alpha\ho 1, 1\ho\beta): (A'\ho B) \oplus(A\ho B') \to A\ho B).\label{tensor2} 
\end{eqnarray}
The extended homology of (\ref{chain}) in dimension 1 will be called {\it the periodic 
product of $\X$ and $\Y$}:
\begin{eqnarray}
\X\ast\Y = \left(\left(
\begin{array}{c}
-1\ho\beta \\
\alpha\ho 1
\end{array}\right)
: A'\ho B'\to Z\right),\end{eqnarray}
where 
\begin{eqnarray}
Z = \ker \left[ \left(
\begin{array}{c}
-1\ho\beta \\
\alpha\ho 1
\end{array}\right)
: (A'\ho B) \oplus 
(A\ho B') \to A\ho B)\right].\label{tensor3}
\end{eqnarray}

It is easy to see that $\X\ho \Y$ and $\X\ast \Y$ are covariant functors of two variables.

\begin{proposition}\label{prop3} 
Let $\ho: \CC\times \CC'\to \CC''$ be a tensor product functor (\ref{tensor}). 
Let $\X\in \ob(\E(\CC))$ and $\Y\in \ob(\E(\CC'))$. Then 
\begin{enumerate}
\item[(a)] $\X\ho \Y$ is projective if both $\X$ and $\Y$ are projective;
\item[(b)] $\X\ast \Y=0$ if $\X$ or $\Y$ is projective;
\item[(c)] $\X\ho \Y$ is torsion if $\X$ or $\Y$ is torsion;
\item[(d)] If $\CC''$ is a finite von Neumann category then $\X\ast \Y$ is torsion for any $\X$ and $\Y$.
\end{enumerate}
\end{proposition}
\begin{proof} Statements (a) and (b) follow directly from the definitions.

Let's prove (c) assuming that $\X=(\alpha: A'\to A)$ is torsion, i.e. $\im \alpha \subset A$ is dense.
From the definition of the tensor product $\ho$ it follows that then the image of
$\alpha\ho 1: A'\ho B\to A\ho B$ is dense and hence from (\ref{tensor2}) we see that $\X\ho\Y$
is torsion. 

It is enough to prove (d), assuming that both $\X$ and $\Y$ are torsion. Let $\X=(\alpha: A'\to A)$
and  $\Y=(\beta: B'\to B)$ with $\alpha$ and $\beta$  injective and having dense images. 
Then $A'$ is isomorphic to $A$ and
$B'$ is isomorphic to $B$ (cf. \cite{F2}, \S 2). Therefore (d) will follow if we can show that 
$Z$ (given by (\ref{tensor3})) is isomorphic to $A\ho B$. The projection of $Z$ on the first coordinate
gives a morphism $Z\to A'\ho B$ which is injective (obviously) and has a dense image (this follows
from Proposition in \S 2 of \cite{F2}). Hence we obtain (using Lemma in \S 2 of \cite{F2})
that $Z$ is isomorphic to $A'\ho B\simeq A\ho B$. 
\end{proof}

\begin{theorem}[Kunneth formula]\label{thm4}
Extended homology $\H_\ast(C\ho C')$ of a tensor product, where $(C, d)$ 
is a chain complex in $\CC$ and $(C', d)$ is a chain complex in $\CC'$, equals
\begin{equation}
\H_n(C\ho C') = \bigoplus_{i+j=n} \H_i(C)\ho\H_j(C') \, \oplus \, 
\bigoplus_{i+j=n-1} \H_i(C)\ast\H_j(C').\label{kunneth}
\end{equation}
\end{theorem}
\begin{proof} Let $Z_i\subset C_i$ and $Z'_i\subset C_i'$ denote the subspaces of cycles. 

We
have the decomposition $C_i=Z_i \oplus Z_i^{\perp}$; the boundary homomorphism vanishes on $Z_i$ and maps $Z_i^{\perp}$ into $Z_{i-1}$. Let's denote by
$D_i$ the short chain complex $D_i=(d: Z_{i+1}^\perp \to Z_i)$, where $Z_i$ stands in degree $i$ and 
$Z_{i+1}^\perp$ stands in degree $i+1$. Then 
$C \simeq \bigoplus_{i=0}^\infty D_i,$
i.e. $C$ is isomorphic to the direct sum of the chain complexes $D_i$.

Similarly, we define chain complexes $D_j'=(d: {Z_{j+1}'}^\perp\to Z_j')$ and 
$C' \simeq \bigoplus_{j=0}^\infty D_j'.$
Hence we obtain
\begin{eqnarray}
C\ho C' \simeq \bigoplus_{i,j} (D_i\ho D_j'), \quad \H_n(C\ho C') = \bigoplus_{i,j} 
\H_n(D_i\ho D_j').\label{sum}
\end{eqnarray}
Now we observe that $D_i$ has nontrivial homology only in dimension $i$ and $\H_i(D_i)=\H_i(C)$;
similarly, $D'_j$ has nontrivial homology only in dimension $j$ and $\H_j(D'_j)=\H_j(C')$. Therefore
$D_i\ho D_j'$ has nontrivial homology only in dimensions $i+j$
and $i+j+1$, and 
\begin{eqnarray}
\H_{i+j}(D_i\ho D_j') = \H_i(C)\ho \H_j(C'), \quad \H_{i+j+1}(D_i\ho D_j') = \H_i(C)\ast \H_j(C')
\label{sum1}\end{eqnarray}
according to our definition of the tensor and periodic products. Formula (\ref{kunneth}) now follows
by combining (\ref{sum}) and (\ref{sum1}).
\end{proof}

\begin{theorem}[Kunneth formula for extended $L^2$ homology]\label{thm5}
Let $X$, $X'$ be finite cell complexes with $\pi=\pi_1(X)$,
$\pi'=\pi_1(X')$. Then
\begin{eqnarray}
\lefteqn{\H_n(X\times X';\ell^2(\pi\times \pi')) \simeq  }\nonumber\\
& & \bigoplus_{i+j=n} \H_i(X;\ell^2(\pi))\ho \H_j(X';\ell^2(\pi'))\oplus \label{kunneth2}\\
&& \bigoplus_{i+j=n-1} \H_i(X;\ell^2(\pi))\ast \H_j(X';\ell^2(\pi')),\nonumber
\end{eqnarray}
where the tensor and periodic products are understood with respect to functor (\ref{tensor1}).
\end{theorem}
\begin{proof} Let $C_\ast(\tilde X)$ and $C_\ast(\tilde X')$ be the cell chain complexes of the universal coverings $\tilde X$ and $\tilde X'$. 
We apply the previous Theorem to chain complexes 
$C=\ell^2(\pi)\tilde\otimes_\pi C_\ast(\tilde X)$ and 
$C'=\ell^2(\pi')\tilde\otimes_{\pi'} C_\ast(\tilde X')$. Note that $C$ is a chain complex in category
$\CC_\pi$ (cf. Example above) and $\H_n(C) = \H_n(X;\ell^2(\pi))$.
Similarly $C'$ is a chain complex in $\CC_{\pi'}$ and $\H_n(C') = \H_n(X';\ell^2(\pi'))$. 
Formula (\ref{kunneth2}) follows from (\ref{kunneth}) using the isomorphism
$\ell^2(\pi)\ho \ell^2(\pi') = \ell^2(\pi\times \pi')$ and the fact that the 
chain complex
$C_\ast(\tilde X)\otimes_{\Z} C_\ast(\tilde X')$ over $\Z[\pi\times \pi']$
is isomorphic to 
$C_\ast(\widetilde {X\times X'})$, where we consider the obvious product cell structure on 
$X\times X'$.
\end{proof}

\bibliographystyle{amsalpha}

\enddocument